\documentclass[12pt]{article}
\usepackage{amssymb,amsmath}

\newtheorem{theorem}{Theorem}
\newtheorem{lemma}{Lemma}

\newtheorem{corollary}{Corollary}

\newcommand{\Z}{{\mathbb Z}}

\newcommand{\eps}{\varepsilon}

\title{$k$-fold sums from a set with few products} 

\author{{\it Dedicated to the memory of Gy\"orgy Elekes} \\ \\ 
Ernie Croot \\ Derrick Hart}

\begin{document}

\maketitle

\section{Introduction}

Before we state our main theorems, we begin with some notation:  given
a finite subset $A$ of some commutative ring, we let $A+A$ denote
the set of sums $a+b$, where $a,b \in A$; and, we let $A.A$ denote
the set of products $ab$, $a,b \in A$.  When three or more sums or
products are used, we let $kA$ denote the $k$-fold sumset 
$A+A+\cdots +A$, and let $A^{(k)}$ denote the $k$-fold product set
$A.A...A$.  Lastly, by $d*A$ we mean the set $A$ dilated by the
scalar $d$, which is just the set $da$, $a \in A$.
\bigskip

The theory of sum-product inequalities has an interesting history, and
began with the theorem of Erd\H os and Szemer\'edi 
\cite{erdos1}, which says that for some $\eps > 0$ and $n > n_0(\eps)$,
we have that for any set $A$ of $n$ real numbers, either the sumset $A+A$
or the product set $A.A$ has at least $n^{1+\eps}$ elements.  
Further improvements to this result were achieved by Nathanson
\cite{nathanson}, Ford \cite{ford}, Elekes \cite{elekes}, and
finally Solymosi \cite{solymosi1} and \cite{solymosi2}.

Another type of theorem that one can prove regarding sums and products
is to assume that either the sumset $A+A$ is near to being as small as
possible (near to $n$), and then to show that $A.A$ must be near to $n^2$;
or, one can suppose that the product set $A.A$ is small, and show that
the sumset $A+A$ is large.  Furthermore, one can consider $k$-fold
sums and products here.  Some quite interesting results along these
lines were produced by Chang \cite{chang2}, \cite{chipeniuk}, 
Elekes-Ruzsa \cite{elekes3}, Elekes-Nathanson-Ruzsa \cite{elekes2}, 
and Jones-Rudnev \cite{jones}.  

There are also some related analogues in finite fields similar to these
just mentioned.  For example, \cite{bourgain2}, \cite{glibichuk}.
\cite{glibichuk2}, \cite{hart} and \cite{vinh}.

Continuing with the characteristic $0$ case, Chang and Bourgain proved
the following results on $k$-fold sums and products: 
Chang \cite{chang1} showed that if $A$ is a set of $n$ integers, and
$|A.A| < n^{1+\eps}$, then the sumset $|kA| \gg_{\eps,k} n^{k-\delta}$,
where $\delta \to 0$ as $\eps \to 0$.  And then Chang and Bourgain
\cite{bourgain} showed that for any $b \geq 1$, there exists 
$k \geq 1$ such that if $A$ is a set of $n$ integers, 
$$
|kA|\cdot |A^{(k)}|\ \gg\ n^b.
$$
In both of these results, we would like to have that they hold for the
real numbers (or even the complex numbers), instead of just the integers.
Unfortunately, this appears to be out of reach at the moment.

The purpose of the present paper is to present some results towards this
end.  Specifically, we will prove the following two theorems.

\begin{theorem} \label{main_theorem1}
For all $h \geq 2$ and $0 < \eps < \eps_0(h)$ we have that the 
following holds for all $n > n_0(h,\eps)$:  if $A$ is a set of 
$n$ real numbers and 
$$
|A.A|\ \leq\ n^{1+\eps},
$$
then
$$
|hA|\ \geq\ n^{\log(h/2)/2\log 2 + 1/2 - f_h(\eps)},
$$
where $f_h(\eps) \to 0$ as $\eps \to 0$.
\end{theorem}

If instead of showing that $kA$ is large, we just want to
show that $k(A.A)$ is large, we can prove a much stronger theorem:

\begin{theorem} \label{main_theorem2} 
Under the same hypotheses on $A$ as in the theorem above, we have
that
$$
|h(A.A)|\ =\ |A.A+A.A+\cdots + A.A|\ >\ n^{\Omega((h/\log h)^{1/3})}.
$$
\end{theorem}

\subsection{Some remarks}

While the methods in the present paper will need substantial modification
to come anywhere near to proving an analogue of \cite{chang1} for the
real numbers, we feel that it might be possible to achieve bounds as
good in Theorem \ref{main_theorem1} above as we have in Theorem
\ref{main_theorem2}.  Although this too will require a lot of work, 
we feel that we have a few good ideas on how to actually achieve it.

It is also worth remarking that we have several different
approaches to proving a theorem of the quality of Theorem 
\ref{main_theorem1}.  In particular, it is possible to use an iterative
argument involving the Szemer\'edi-Trotter theorem, a Szemer\'edi cube
lemma similar to Lemma \ref{box_lemma} below, and some ``energy arguments'',
to achieve similar such bounds.  However, it is not as easy to see how one
might go about modifying such an ``incidence proof'' of 
Theorem \ref{main_theorem1} to achieve bounds as good as in 
Theorem \ref{main_theorem2}.

\section{Preliminary lemmas and results} \label{prelim_section}

First, in the proofs of both theorems, we will assume that all the
elements of $A$ are positive.  The reason we can assume this is that
either at least $(n-1)/2$ elements of $A$ are all positive, or at least
$(n-1)/2$ are all negative.  If we are in the negative case here, we
just let $A'$ be the negative of these negative elements; and otherwise, 
we just let $A'$ be these $(n-1)/2$ positive elements.  Then, we simply prove
our theorem using $A'$ in place of $A$.  The effect of the lost factor
of $2$ will be negligible.  
\bigskip

The proof of Theorem \ref{main_theorem2} will require the following
result of Wooley \cite{wooley} (see also Borwein-Erd\'elyi-K\'os 
\cite{borwein} for some related results).

\begin{theorem} \label{wooley_theorem}
For every $k \geq 1$, there exist two distinct sets 
$$
\{x_1,...,x_s\},\ \{y_1,...,y_s\}\ \subseteq\ \Z
$$
of 
$$
s\ <\ (k^2/2)(\log k + \log\log k + O(1))
$$
(when $k=1$ we just delete the $\log\log k$ term)
distinct integers in each, such that
$$
\sum_{i=1}^s x_i^j\ =\ \sum_{i=1}^s y_i^j,\ {\rm for\ all\ } j=1,...,k,
$$
but that
$$
\sum_{i=1}^s x_i^{k+1}\ \neq\ \sum_{i=1}^s y_i^{k+1}.
$$
\end{theorem}

For the purposes of our paper, we actually require the following corollary
of this theorem.

\begin{corollary} \label{wooley_corollary}
For all integers $j \geq 1$, there exists a monic polynomial $f(x)$, having
only the coefficients $0,1,$ and $-1$, having at most 
$$
(j^2/2)(\log j + \log\log j + O(1))
$$
(again, if $j=1$ we just delete the $\log\log j$ term)
non-zero terms, such that $f(x)$ vanishes at $x=1$ to order $j$, but
not to order $j+1$. 

Since the number of terms of this polynomial depends only on $j$, it
follows that if one performs a Taylor expansion of this polynomial
about $x=1$, one will find that
$$
f(x)\ =\ \sum_{i=j}^d c_i (x-1)^i,\ d={\rm deg}(f),
$$ 
where $c_j \neq 0$, and each $c_i$ in turn is either $0$ or its absolute 
value can be bounded from below by some function of $j$ alone 
(the degree $d$ depends on $j$).
\end{corollary}

\noindent {\bf Proof of the Corollary.}  Basically, we just use the 
well-known fact that a polynomial
$$
f(x)\ =\ \sum_{i=1}^s x^{x_i}\ - \sum_{i=1}^s x^{y_i}.
$$
vanishes to order $j$ at $x=1$ if and only if $f$ and its first 
$j-1$ derivatives vanish at $x=1$, where the $\ell$th derivative 
evaluated at $1$ looks like
$$
\sum_{i=1}^s x_i(x_i-1)\cdots (x_i-\ell+1)\ -\ \sum_{i=1}^s y_i (y_i-1)
\cdots (y_i - \ell+1).
$$
Clearly, having all these be $0$, for $\ell=0,1,...,j-1$, is equivalent
to having a solution to the ``Tarry-Escott Problem'' considered by 
Wooley in Theorem \ref{wooley_theorem} above.  

One small remaining point to consider is the fact that some of the
$x_i$'s and $y_i$'s could be negative (meaning that the $f$ above is
a Laurent polynomial, not a polynomial).  That is easily fixed by
multiplying $f$ by an appropriate power of $x$, which does not affect
the vanishing properties at $x=1$.
\hfill $\blacksquare$
\bigskip

Another major theorem that we will require is the Ruzsa-Plunnecke
inequality \cite{ruzsa}.

\begin{theorem}  Suppose that $A$ is a finite subset of an additive
abelian group.  Then, if 
$$
|A+A|\ \leq\ K|A|,
$$
we will have that 
$$
|kA - \ell A|\ =\ |A+A+\cdots + A - A - \cdots - A|\ \leq\ K^{k+\ell} |A|.
$$
\end{theorem}

We will also require the following basic lemmas.

\begin{lemma} \label{already_good_lemma}  
Suppose that $A$ is a set of $m^2$ positive real numbers, say they are
$$
0\ <\ a_1\ <\ \cdots\ <\ a_{m^2},
$$
such that no dyadic interval $[x,2x]$ contains $m$ or more of 
the $a_j$'s.  Then,
$$
|kA|\ \gg_k\ m^k.
$$
\end{lemma}

\noindent {\bf Proof of the lemma.}
Let 
$$
B\ :=\ \{a_1, a_{2m+1}, a_{4m+1}, ..., a_{m^2-m+1}\}.
$$
We claim that all the sums
$$
b_1 + \cdots + b_k,\ b_i \in B,\ b_1 < \cdots < b_k,
$$
are distinct, which would prove the lemma.

To see this, suppose we had
$$
b_1 + \cdots + b_k\ =\ b'_1 + \cdots + b'_k,
$$
and suppose without loss that $b_k \leq b'_k$.  If $b_k < b'_k$,
then $b_k < b'_k/2$, and we have
$$
b'_k\ =\ b_1 + \cdots + b_k - b'_1 - \cdots - b'_{k-1}
\ \leq\ b_1 + \cdots + b_k\ <\ 2b_k\ <\ b'_k, 
$$
contradiction.  So, we can delete $b_k$ and $b'_k$ from both sides; 
and then, repeating the argument, we get 
$b_i = b'_i$, $i=1,2,...,k$, and we are done.
\hfill $\blacksquare$
\bigskip

The following lemma is basically a generalization of a result in
\cite{croot}.

\begin{lemma} \label{translate_lemma} 
For every $r, \ell \geq 1$, $0 < c_1 < c_1(r,\ell)$, 
$0 \leq c_2 < c_2(r,\ell,c_1)$ and 
$0 < \eps < \eps(r,\ell,c_1,c_2)$, the following holds for all 
$n$ sufficiently large:  Suppose $A$ is a set of $n$ real numbers satisfying
$$
|A.A|\ <\ n^{1+\eps},
$$
and suppose that 
$$
A'\ \subseteq\ A,\ B\ \subseteq\ A^{(r)}/A^{(r)},
$$
satisfy 
$$
|B|\ >\ n^{c_1},\ {\rm and\ } |A'|\ >\ n^{1-c_2}.
$$
Then, there are 
$$
n^{2-O(c_2 \ell + r\ell^2 \eps)}\ {\rm pairs\ } (a_1,a_2) \in A' \times A', 
$$
such that if we let 
$$
x\ =\ a_1 a_2,
$$
then there exist
$$
n^{-O(c_2 \ell + r \ell^2 \eps)} |B|^2 
$$
pairs 
$$
(b_1,b_2)\ \in\ B \times B,\ b_1 > b_2,
$$
such that if we let
$$
d\ =\ b_1/b_2,
$$
then
$$
x,\ dx,\ d^2x,\ ...,\ d^\ell x\ \in\ A'.A'.
$$
\end{lemma}

\noindent {\bf Proof of the lemma.}  First, by Ruzsa-Plunnecke, we have that
$$
|A'/A'|\ \leq\ |A/A|\ <\ n^{1+2\eps}.
$$
And then, by a simple pigeonhole argument, we have that there exists 
$t_1,t_2 \in A'$ such that if 
$$
A''\ :=\ (A'/t_1) \cap (t_2/A'),
$$
then
$$
|A''|\ >\ n^{1- O(c_2)},
$$
for $c_2 > 0$ small enough in terms of $\eps$.
Of course, this means that if we let $t = t_1t_2$ then
$$
A'''\ :=\ A' \cap (t/A')\ {\rm satisfies\ } |A'''|\ =\ |A''|\ >\ n^{1-O(c_2)}.
$$
Note that
$$
A'''/A'''\ \subseteq\ t^{-1}*(A'.A').
$$
\bigskip

What we will show is that for lots of pairs 
$$
(e_1,e_2)\ \in\ A''' \times A''',\ {\rm and\ }y = e_1/e_2,
$$
the set $A'''/A'''$ contains ``lots'' of progressions
\begin{equation} \label{yd_equation}
y,\ dy,\ d^2 y, ...,\ d^\ell y,\ {\rm where\ } d \in B/B.
\end{equation}
If we can do this, then 
$$
x,\ dx,\ d^2 x,\ ...,\ d^\ell x\ \in\ A'.A',\ {\rm where\ } x = ty.
$$
\bigskip

Showing that $A'''/A'''$ contains such a progression amounts to showing 
that the following system has non-trivial solutions:
\begin{equation} \label{the_good_system}
{b_1 \over b_2} {e_1 \over e_2}\ =\ {e_3 \over e_4},\ 
{b_1^2 \over b_2^2} {e_1 \over e_2}\ =\ {e_5 \over e_6},\ ...,\ 
{b_1^\ell \over b_2^\ell} {e_1 \over e_2}\ =\ {e_{2\ell+1} \over 
e_{2\ell +2}}.
\end{equation}
Another way to write this is
\begin{eqnarray*}
b_1 e_1 e_4\ &=&\ b_2 e_2 e_3 \\
b_1^2 e_1 e_6\ &=&\ b_2^2 e_2 e_5 \\
&\vdots& \\
b_1^\ell e_1 e_{2\ell+2}\ &=&\ b_2^\ell e_2 e_{2\ell+1}.
\end{eqnarray*}

Note that both sides of the equations belong to 
$$
t^{-2} *(B^{(\ell)} (A')^{(4)})\ \subseteq\ t^{-2}*(A')^{(r\ell+4)},
$$
which has at most
$$
n^{1 + O(r\ell \eps)}
$$
elements, by Ruzsa-Plunnecke.  So, there are 
$$
\leq\ n^{\ell + O(r \ell^2 \eps)}
$$
values that the $\ell$ expressions on the left-hand-side in this
system, in total, can take on, and the same goes for the right-hand-side.

Now, since there are at least 
$$
n^{\ell + 1 - O(c_2 \ell)} |B|
$$
choices for
$$
e_1,\ e_4,\ e_6,\ e_8,\ ...,\ e_{2\ell+2},\ b_1,
$$
that make up the left-hand-sides -- similarly, the right-hand-sides -- 
it is clear that we get at least 
$$
n^{\ell + 2 - O(c_2 \ell + r \ell^2 \eps)} |B|^2
$$
solutions.  For $c_1,c_2$ small enough in terms of $r, \ell$ and for 
$\eps$ small enough in terms of $r,\ell,c_1,c_2$, this exceeds the
number of tirival solutions, which is 
$$
|B| \sum_{e \in A'''/A'''} r(e)^\ell\ \leq\ |B|\cdot |A'''/A'''|\cdot 
|A'''|^\ell\ =\ n^{\ell + 1 + O(\ell c_2)}|B|,
$$
where $r(e)$ is the number of representations $e = e'/e''$, where
$e',e'' \in A'''$.

So, there are at least 
$$
n^{\ell + 2 - O(c_2 \ell + r \ell^2 \eps)} |B|^2
$$
non-trivial solutions, when $c_1,c_2$ and $\eps > 0$ are small enough.
So, for the average tuple 
$$
(e_4, e_6, ..., e_{2\ell+2})\ \in\ (A''')^{\ell},
$$
we have that there are at least
$$
n^{2-O(c_2 \ell + r \ell^2 \eps)} |B|^2
$$
four-tuples
$$
(e_1,e_2,b_1,b_2)\ \in\ A''' \times A''' \times B \times B,
$$
such that the system has a solution (which must be unique, since
the remaining $c_i$'s are determined exactly).  Clearly this proves the
lemma.  Note that guaranteeing that we can work with $b_1 > b_2$
can be guaranteed simply by taking reciprocals in (\ref{the_good_system}).
\hfill $\blacksquare$
\bigskip

And now we state two more general-purpose lemmas, the first of which
is perhaps better known as a the ``Szemeredi cube lemma'', 
and is used in the proof of Theorem \ref{main_theorem1} only, 
while the other lemma is used in the proofs of both theorems.

\begin{lemma} \label{box_lemma}  The following holds for all $k \geq 2$,
$0 < c < c_0(k)$, $0 < \eps < \eps_0(k,c)$ and $n > n_0(k, c, \eps)$:
suppose that $A$ is a set of $n$ real numbers such that 
$$
|A.A|\ <\ n^{1+\eps},
$$
and suppose that $B \subseteq A$ satisfies
$$
|B|\ \geq\ n^c.
$$
Then, there exists 
$$
\theta_1,\ ...,\ \theta_k\ \in\ B/B,\ {\rm each\ }\theta_i > 1,
$$
such that for at least 
$$
n^{1-O_k(c)}
$$
values $d \in A$ we have that all the numbers
$$
\theta_1^{\gamma_1} \cdots \theta_k^{\gamma_i} d,\ {\rm where\ each\ }
\gamma_i\ \in\ \{0,1\},
$$
belong to $A$.  (Note that for each such $d$, this means that 
$2^k$ different elements belong to the set $A$.)
\end{lemma}

\noindent {\bf Proof of the lemma.}  The proof is inductive:  we will
construct
$$
D_0\ :=\ A,\ D_1,\ ...,\ D_k\ \subseteq\ A,
$$
such that 
$$
D_i\ =\ D_{i-1} \cap (\theta_i^{-1} * D_{i-1}),\ i=1,2,...,k,
$$
where $\theta_i \in B/B$, $\theta_i > 1$,
is chosen greedily to maximize $D_i$, given $D_{i-1}$.  

Suppose that we have already shown that
$$
|D_{i-1}|\ =\ n^{1-O_k(c)}.
$$
Then, consider the product set 
$$
BD_{i-1}\ =\ \{b d\ :\ b \in B,\ d \in D_{i-1}\}\ \subseteq\ A.A.
$$
Since $|A.A| < n^{1+\eps}$, we have that 
$$
|BD_{i-1}|\ \leq\ n^{1+\eps},
$$
and therefore we easily see that there exists 
$$
s,t \in B,\ s \neq t,\ s/t \neq \theta_1,...,\theta_{i-1},\ s > t,
$$
such that 
$$
|(s*D_{i-1}) \cap (t*D_{i-1})|\ >\ n^{1-O_k(c)},
$$
for $\eps > 0$ sufficiently small in terms of $k,c$.

So, letting $\theta_i = s/t > 1$, we are done, because
$$
|D_{i-1} \cap (\theta_i^{-1}*D_{i-1})|\ =\ |(s*D_{i-1}) \cap (t*D_{i-1})|
\ >\ n^{1-O_k(c)},
$$
as claimed.
\hfill $\blacksquare$
\bigskip

\begin{lemma} \label{distinct_sums_lemma} 
Suppose that $C$ is a set of real numbers, and
$$
1\ =\ \delta_0\ >\ \delta_1\ >\ \delta_2\ >\ \cdots\ >\ \delta_{k-1}\ >\ 0
$$
are positive real numbers such that if we define the ratios
$$
\alpha_i\ :=\ \delta_i/\delta_{i-1},\ i=1,2,...,k-1,
$$
then for all pairs 
$$
c,d\ \in\ C,\ c > d,
$$
we have
$$
c/d - 1\ >\ 2k \alpha_2,\ ...,\ 2k \alpha_{k-1}.
$$
Next, partition $C$ into any disjoint sets
$$
C\ =\ C_1\ \cup\ C_2\ \cdots\ \cup\ C_k,
$$
where for $i < j$ we have that every element of $C_i$ is greater than
every element of $C_j$.  Let us express this as
$$
C_i\ >\ C_j,\ {\rm for\ }i < j.
$$  

Then, we have that all sums
$$
c_1 + c_2 \delta_1 + \cdots + c_k \delta_{k-1},\ c_1,...,c_k \in C,
$$
are distinct.
\end{lemma}

\noindent {\bf Proof of the lemma.}  Suppose that, on the contrary, two of
these sums are equal.  Then, it would mean that 
\begin{equation} \label{both_sides}
c_1 + c_2 \delta_1 + \cdots + c_k \delta_{k-1}\ =\ 
c'_1 + c'_2 \delta_1 + \cdots + c'_k \delta_{k-1}.
\end{equation}
Suppose without loss that $c_1 \geq c'_1$.  Now let us suppose
that, in fact, $c_1 > c'_1$.  Then, we have that
$$
c_1/c'_1 - 1\ =\ \sum_{i=2}^k (c'_i/c'_1 - c_i/c'_1) \delta_{i-1}. 
$$
From the fact that $C_1 > C_2,...,C_k$, we have that the right-hand-side 
here is bounded from above in absolute value by
$$
\sum_{i=2}^k 2 \delta_i\ \leq\ 2k\delta_2\ \leq\ 2k \alpha_2\ <\ c_1/c'_1 - 1,
$$
which is impossible.  We conclude that $c_1 = c'_1$.

Now suppose for proof by induction we have shown that 
$$
c_i = c'_i,\ i=1,2,3,...,j,\ {\rm where\ } j \leq k-1.
$$
We now show that
$$
c_{j+1}\ =\ c'_{j+1},
$$
which would clearly prove the lemma (we get $c_k = c'_k$ for free
once the other $c_i = c'_i$ are established).

We begin by deleting the terms $c_i \delta_{i-1}$ and
$c'_i \delta_{i-1}$ from both sides of (\ref{both_sides}), for
$i=1,2,...,j$.  So, we are left with
$$
\sum_{i=j+1}^k c_i \delta_{i-1}\ =\ \sum_{i=j+1}^k c'_i \delta_{i-1},
$$
which can be rewritten as
\begin{equation} \label{both_sides2}
(c_{j+1}/c'_{j+1} - 1)\ =\ \sum_{i=j+2}^k (c'_i/c'_{j+1} - c_i/c'_{j+1})
\delta'_{i-1},
\end{equation}
where
$$
\delta'_i\ :=\ \delta_i/\delta_j\ \leq\ \alpha_i,\ i=j+1,...,k-1. 
$$

We assume without loss that $c_{j+1} \geq c'_{j+1}$.
If, in fact, $c_{j+1} > c'_{j+1}$, then the absolute value
of the right-hand-side of (\ref{both_sides2}) is 
clearly bounded from above by 
$$
2k \alpha_{j+1}\ <\ c_{j+1}/c'_{j+1} - 1,
$$
which is a contradiction.  We conclude that $c_{j+1} = c'_{j+1}$,
and therefore the induction step is proved, as is the lemma.
\hfill $\blacksquare$

\section{Proof of Theorem \ref{main_theorem1}}

We suppose that the elements of $A$ are
$$
a_1\ <\ a_2\ <\ \cdots\ <\ a_n,
$$
which we assume are all positive by the remarks at the beginning of 
section \ref{prelim_section}.
 
Let $0 < \delta < 1/2$ be some parameter that we will choose as small as needed
later, and let $k \geq 2$ be some parameter that we will let depend
on $h$ later.  Let 
$$
s\ :=\ \lfloor n^\delta \rfloor,
$$
and set
$$
B\ :=\ \{a_j, a_{j+1}, ..., a_{j+s}\},
$$
where $j$ is chosen so that 
\begin{equation} \label{minimal_gap}
a_{j+s}/a_j\ {\rm is\ minimal}.
\end{equation}
We may assume that $B$ lies in some interval
$$
B\ \subseteq\ [x,2x],
$$
since otherwise each consecutive block of $s+1$ elements 
of $A$ lies in its own interval of this form, disjoint from those corresponding
to other blocks of elements of $A$; and therefore, by Lemma 
\ref{already_good_lemma}, we could conclude that 
$$
|kA|\ \geq\ n^{\Omega(\sqrt{k})},
$$
which would prove our theorem.

Having $B$ lie in a dyadic interval implies that all ratios $\theta = b_2/b_1$,
$b_1,b_2 \in B$, satisfy
$$
\theta - 1\ \in\ [0,1).
$$ 
\bigskip

Next, we let $c = \delta$ and apply Lemma \ref{box_lemma}, and let
$\theta_1, ..., \theta_{k-1} \in B/B$ denote the numbers that 
result from this lemma (using $k-1$ in place of $k$); and, let
$C_0$ denote the set of all $n^{1 - O_k(\delta)}$ elements $d \in A$ 
that the lemma produces, and write
$$
C_0\ :=\ \{c_1, ..., c_{n'}\},\ c_1 < c_2 < \cdots < c_{n'},\ 
n'\ >\ n^{1-O_k(\delta)}.
$$  
Then, let 
$$
C\ :=\ \{c_1, c_{1 + sk2^k}, c_{1 + 2 sk2^k}, c_{1 + 3sk2^k},\ ...\}
$$
be the set of every $sk2^k$th element of $C_0$.  Note that since
the elements $\theta_i \in B/B$, and $B$ satisfies (\ref{minimal_gap}),
we have that
$$
\theta_j^{2k}\ <\ c_2/c_1,\ c_1,c_2 \in C,\ c_2 > c_1,
$$
and therefore
$$
2k(\theta_j - 1)\ <\ \theta_j^{2k} - 1\ <\ c_2/c_1 - 1.
$$
It follows that if we let
$$
\delta_i\ :=\ (\theta_1 - 1)(\theta_2 - 1) \cdots (\theta_i - 1),\ 
i=1,2,...,k-1,
$$ 
then we have for $i=1,2,...,k-2$ that
$$
\delta_{i+1}/\delta_i\ =\ \theta_{i+1} - 1\ <\ (c_2/c_1 - 1)/2k.
$$

We almost are ready to apply Lemma \ref{distinct_sums_lemma} -- all we
have to do is partition $C$, which we do simply by letting 
$$
C\ =\ C_1 \cup \cdots \cup C_k,
$$
where $C_1$ consists of the largest $\lfloor |C|/k\rfloor$ elements of 
$C$, $C_2$ consists of the next largest $\lfloor |C|/k\rfloor$ elements of 
$C$, and so on. 

Lemma \ref{distinct_sums_lemma} now tells us that all the sums
\begin{equation} \label{d_sums}
c_1 + c_2 \delta_1 + \cdots + c_k \delta_{k-1},\ c_i \in C_i,
\end{equation}
are distinct.  This then results in 
$$
\gg\ (|C|/k)^k\ \gg_k\ n^{k(1-O_k(\delta))}.
$$
distinct sums.  

These sums, in turn, can be re-written as just sums and differences 
of elements from $A$ as follows:  by expressing the $\delta_{i-1}$ back in 
terms of the $\theta_j$'s, we find that 
$$
c_i \delta_{i-1}\ =\ c_i (\theta_1 - 1)\cdots (\theta_{i-1}-1)
\ =\ (-1)^{i-1} c_i + (-1)^{i-2} \theta_1 c_i + \cdots
$$ 
Each term here looks like
$$
\pm c_i \theta_1^{\gamma_1} \cdots \theta_{i-1}^{\gamma_{i-1}},\ {\rm where\ }
\gamma_j \in \{0,1\},
$$
and we know from our use of Lemma \ref{box_lemma} that all such numbers
belong to $\pm A$.  

It is easy to see, then, that all the sums (\ref{d_sums}) can be 
re-expressed as subsets of $KA - LA$, where
$$
K,L\ <\ 2^{k-2} + 2^{k-3} + \cdots + 1\ <\ 2^{k-1}. 
$$
So,
$$
|2^{k-1}A|^2\ >\ |KA|\cdot |LA|\ \geq\ |KA - LA|\ >\ n^{k(1 - O_k(\delta))},
$$
from which it follows that upon letting $h = 2^{k-1}$,
$$
|hA|\ \geq\ n^{\log(h)/2\log 2 + 1/2 - g_h(\delta)},
$$
where $g_h(\delta) \to 0$ as $\delta \to 0$.  Of course, this only 
works for when $h$ is a power of $2$; by bounding general $h$ between
two consecutive powers of $2$, we can conclude that
$$
|hA|\ \geq\ n^{\log(h/2)/2\log 2 + 1/2 - g_h(\delta)}.
$$

This completes the proof of our theorem, by choosing $\delta > 0$ small
enough, and then choosing $\eps > 0$ even smaller as needed.

\section{Proof of Theorem \ref{main_theorem2}} 

Write out the elements of $A$ in incresing order as
$$
a_1\ <\ a_2\ <\ \cdots\ <\ a_n,
$$
which we assume are all positive by the remarks at the beginning
of section \ref{prelim_section}.   

Let $\delta > 0$ be some parameter that we will 
choose later as function of $h$ alone, and let $s$ and $B$ be
as in the beginning of the proof of Theorem \ref{main_theorem1}.
In fact, we may assume that 
$$
B\ \subseteq\ [x,(1 + 1/\gamma(h)) x],
$$
for any function $\gamma(h) > 0$ we please.  The reason is that 
if this minimal set $B$ lies in this interval, and if $a_j$ is 
the smallest element of $B$, then 
$$
[a_j,\ a_{j + L s}]\ \supseteq\ [x, (1 + \gamma(h))^L x]\ \supseteq\ 
[x,2x],
$$
for $L$ large enough.  And so, we may again deduce (as in the
proof of Theorem \ref{main_theorem1}), using Lemma \ref{already_good_lemma}, 
that
$$
|hA|\ \gg\ n^{\Omega(h^{1/2})}, 
$$
which would prove our lemma.
\bigskip

As in the proof of Theorem \ref{main_theorem1}, we may assume that
the elements of $B$ lie in some dyadic interval $[x,2x]$, upon applying
Lemma \ref{already_good_lemma}.

Let $k \geq 2$ be some parameter that is to depend on $h$, that we 
will choose later.  

We now apply Lemma \ref{translate_lemma} using $\ell = M$, $r=1$, 
$c_1 = \delta$, $A' = A$ (so $c_2 = 0$) and $\eps > 0$ as small as 
needed in terms of $M$ and $\delta$, where the precise value of $M$ will be 
determined below, and will depend only on $h$.  So, there exists 
$$
\theta\ \in\ B/B,\ \theta\ >\ 1,
$$
such that for at least 
$$
|A|^{2 - O(\ell^2 \delta)}
$$
pairs
$$
(a_1,a_2)\ \in\ A,
$$
we have that if we let $y = a_1a_2$, then
\begin{equation} \label{xtheta}
y,\ y \theta,\ ...,\ y \theta^\ell\ \in\ A.A.
\end{equation}
So, there exists $a_1 \in A$ and $\theta \in B/B$, $\theta > 1$, 
such that there are
\begin{equation} \label{A'count}
|A|^{1 - O(\ell^2 \delta)}
\end{equation}
values $a_2 \in A$ such that for $y=a_1a_2$ we have that 
(\ref{xtheta}) holds.  

If we let the special elements $a_2 \in A$ be
$$
\{a_{h_1}, ..., a_{h_{n'}}\},\ n' > n^{1 - O(\ell^2 \delta)},
$$
then we define 
\begin{equation} \label{A''def}
C\ :=\ \{a_{h_1}, a_{h_{\lfloor \sqrt{n} \rfloor +1}}, 
a_{h_{2\lfloor \sqrt{n} \rfloor + 1}}, ... \}\ \subseteq\ 
\{a_{h_1},...,a_{h_{n'}}\}.
\end{equation}
Note that $C$ is basically a ``well-separated'' subset of those special
elements $a_2\in A$; and, in fact, if 
$\delta < 1/2$, so that $s > n^{1/2}$, 
they are so well-separated that that all ratios 
$c_2/c_1$, $c_1,c_2 \in C$, $c_2 > c_1$, have the property that  
\begin{equation} \label{cctheta}
c_2/c_1\ >\ \theta.
\end{equation}
Later, we will prove and make use of an even stronger such inequality. 
\bigskip

Now we apply Corollary \ref{wooley_corollary}, letting 
$$
f_1(x),\ ..., f_{k-1}(x)
$$
be polynomials having at most $(j^2/2)(\log j + \log\log j + O(1))$,
$j=1,...,k-1$, terms each, each with coefficients only
$0,1,$ or $-1$, that vanish at $x=1$ to the orders 
$1,2,3,...,k-1$, respectively.  Then we let 
$$
M\ =\ \max({\rm deg}(f_1),...,{\rm deg}(f_{k-1})),
$$
which is a parameter that came up earlier in the proof of the present
theorem (Theorem \ref{main_theorem2}).  Note that $M$ does not 
depend on $n$ -- it depends on $k$, and therefore on $h$.

Next, we set
$$
\delta_i\ :=\ f_i(\theta),\ i=1,2,...,k-1.
$$
Since $f_{i+1}(x)/f_i(x)$ vanishes at $x=1$ to order $1$, and since 
$\theta \in B/B$ and $\theta \in [1,2)$, we have that if we set 
$$
\alpha_i\ :=\ \delta_i/\delta_{i-1},\ i=2,3,...,k-1,
$$
then for every $c_1,c_2 \in C$, $c_2 > c_1$, we have from 
(\ref{cctheta}) that 
$$
c_2/c_1 - 1\ >\ \theta-1\ \gg_k\ \alpha_i\ >\ 0.
$$
(Note that the implied constant here depends on the sizes of the 
coefficients $c_i$ in Corollary \ref{wooley_corollary}, and we
know that these coefficients are rational numbers that depend on $k$.)
But, in fact, if $\delta > 0$ is small enough, then for 
$c_1,c_2 \in C$, $c_2 > c_1$, we can
assume that for any function $\gamma(k)$ of $k$, 
$$
\theta^{\gamma(k)}\ <\ c_2/c_1;
$$
and so, we may assume
$$
c_2/c_1 - 1\ >\ 2k \alpha_i\ >\ 0.
$$
So, if we let $C_1$ be the largest $\lfloor |C|/k \rfloor$ elements of $C$, 
$C_2$ be the second largest $\lfloor |C|/k\rfloor$ elements of $C$,
and so on, down to $C_k$, then upon applying Lemma \ref{distinct_sums_lemma},
we have that all sums 
\begin{equation} \label{poly_C_sum}
a_1 + a_2 f_1(\theta) + \cdots + a_k f_k(\theta),\ a_i \in C_i,
\end{equation}
are distinct.  Since each $C_i$ satisfies
$$
|C_i|\ \gg_k\ n^{1/3},
$$
for $\delta > 0$ small enough (and $\eps > 0$ small enough in terms
of $\delta$ and $h$), we deduce that this produces $n^{\Omega(k)}$ 
distinct sums.  Now, because 
$$
a_i, a_i \theta,\ ...,\ a_i \theta^M\ \in\ A.A,
$$
by design, we have that upon expanding out these polynomials $f_i$
in (\ref{poly_C_sum}) into powers of $\theta$, we find that these
$n^{\Omega(k)}$ sums are, in fact, subsets of $KA - L A$, where
$$
K,L\ \leq\ (1^2 + 2^2 + \cdots + (k-1)^2/2)(\log k + \log\log k + O(1))\ 
\ll\ k^3 \log k.
$$
It follows that 
$$
|(ck^3 \log k) A|^2\ \geq\ |KA - LA|\ \geq\ n^{\Omega(k)}. 
$$
This clearly proves the theorem upon letting $k \gg (h/\log h)^{1/3}$.

\section{Acknowledgements}

We would like to thank P. Borwein and T. Erd\'elyi for useful discussions
about their paper with G. K\'os listed below.

\end{document}